\documentclass{article}
\usepackage{amssymb,amsmath,amsthm}
\newtheorem{theorem}{Theorem}

\newtheorem{lemma}{Lemma}

\newtheorem{remark}{Remark}
\date{}
\numberwithin{equation}{section} \numberwithin{theorem}{section}
\numberwithin{lemma}{section} \numberwithin{corollary}{section}
\numberwithin{remark}{section} \numberwithin{proposition}{section}
\numberwithin{definition}{section}
\begin{document}
\newcommand{\n}{\noindent}
\newcommand{\vs}{\vskip}

\title{ Porosity  of the Free Boundary in the Singular p-Parabolic Obstacle Problem}

\author{Abdeslem Lyaghfouri\\
American University of Ras Al Khaimah\\
Ras Al Khaimah, UAE} \maketitle

\begin{abstract}
In this paper we establish the exact growth of the solution of
the singular quasilinear p-parabolic obstacle problem near the free  
boundary from which we deduce its porosity.
\end{abstract}

\begin{flushleft}
\end{flushleft}

\begin{flushleft}
\end{flushleft}

\section{Introduction}\label{S:intro}

\n Let  $\Omega$ be an open bounded domain of
$\mathbb{R}^n$, $n\geqslant 2$, $T>0$.
We consider the following problem
\begin{equation*}(P)
\begin{cases}
& \text{ Find }\,\,u\in L^p(0,T; W^{1,p}(\Omega))
\text{ such that } :\\
& (i) \quad  u\geqslant 0\quad\text{in }\Omega_T=\Omega\times(0,T),\\
& (ii)\quad L_p(u)=u_t-\Delta_p u=-f(x)
      \qquad \text{in}\quad \{u>0\},\\
& (iii)\quad u=g \quad  \text{on }\quad \partial_p\Omega_T=(\Omega\times\{0\})\cup(\partial\Omega\times(0,T)),
\end{cases}
\end{equation*}

\n where $p>1$, $\Delta_p$ is the $p$-Laplacian defined by $\Delta_p u=div\big(|\nabla u|^{p-2}\nabla u\big)$,
and $f$, $g$ are functions defined in $\Omega_T$ and satisfying for 
two positive constants $\lambda_0$ and  $\Lambda_0 $ 
\begin{eqnarray}\label{e-1.1}
0< \lambda_0 \leqslant f \leqslant \Lambda_0 \quad \hbox{ a.e. in}~~\Omega_T.
\end{eqnarray}
Moreover we assume that 
\begin{eqnarray}\label{e-1.2-4}
&&f~~\text{ is non-increasing}~~\text{ in}~~t.\\
&&g(x,0)=0~~ \hbox{ a.e. in}~~\Omega.\\
&&g~~\text{ is non-decreasing}~~\text{ in}~~t.
\end{eqnarray}

\vs 0.2cm\n The variational formulation of the problem $(P)$ is given by
\begin{equation*}(VP)
\begin{cases}
& \text{ Find }\,\,u\in K_g=\{v\in V^{1,p}(\Omega_T)~/~v=g\text{ on } \partial_p\Omega_T,~~v\geqslant0 \text{ a.e. in } \Omega_T~
\}\\
&\text{ such that for all } h>0 \text{ and } t<T-h:\\
& \displaystyle{\int_\Omega \partial_t u_h(v-u)dx+\int_\Omega \big(|\nabla u|^{p-2}\nabla u\big)_h.\nabla(v-u)dx+\int_\Omega f_h (v-u)dx\geqslant 0},\\
& \text{a.e. in } t\in(0,T), \text{ and for all } v\in K_g,
\end{cases}
\end{equation*}
where $$V^{1,p}(\Omega_T)=L^\infty(0,T; L^1(\Omega))\cap L^p(0,T; W^{1,p}(\Omega)),$$
and $v_h$ is the Steklov average of a function $v$ defined by
\begin{eqnarray*}
&&v_h(x,t)={1\over h}\int_t^{t+h}v(x,s)ds,\quad\text{ if }  t\in(0,T-h]\\
&&v_h(x,t)=0,\quad\text{ if }  t>T-h.
\end{eqnarray*}

\vs 0.2cm\n Let us recall the following existence and uniqueness theorem of the 
solution of the problem $(VP)$ \cite{[S]}.
\vs 0,2cm \n
\begin{theorem}\label{t1.1} Assume that $f$ and $g$ satisfy (1.1)-(1.4). Then there
exists a unique solution $u$ of the problem $(VP)$ which satisfies
\begin{align*}
   &  0\leqslant u \leqslant M=\| g\|_{\infty,\Omega_T} \quad \hbox{in }\Omega_T.\\
   &  u_t \geqslant 0\quad \hbox{in } \{u>0\}.\\
   &  f \chi_{\{u>0\}} \,\,\leqslant \,\, \Delta_p u-u_t\leqslant f \quad \hbox{a.e.  in }\Omega_T.
\end{align*}
 \end{theorem}

\vs 0,2cm \n
\begin{remark}\label{r1.1} 
We deduce from (1.5) and (1.7) \cite{[Ch]} that we have $u\in C_{loc}^{0,\alpha}(\Omega_T)\cap C_{x,loc}^{1,\alpha}(\Omega_T)$ 
for some $\alpha\in(0,1)$.
 \end{remark} 
    
\vs 0.2cm\n The main result of this paper is the next theorem.
\vs 0.2cm \n
\begin{theorem}\label{t1.2} Assume that $1<p<2$ and that $f$ and $g$ satisfy (1.1)-(1.4),
and let $u$ be a solution of $(VP)$.  Then for every compact set $K\subset\Omega_T$, 
the intersection  $(\partial\{u>0\})\cap K\cap\{t=t_0\} $ is porous in $\mathbb{R}^{n}$ 
with porosity constant depending only on  $n$, $p$, $\lambda_0$, $\Lambda_0$,
$dist(K,\partial_p \Omega_T)$, and $\| g\|_{\infty,\Omega_T}$.
\end{theorem}

\vs 0,2cm\n We recall that a set $E\subset \mathbb{R}^n$ is
called  porous  with porosity $\delta$, if there is an $r_0>0$
such that
$$ \forall x\in E , \quad \forall r\in (0,r_0),\quad \exists y\in
\mathbb{R}^n\quad \hbox{ such that } \quad B_{\delta r}(y)\subset
B_{ r}(x)\setminus E.$$

\n A porous set has  Hausdorff dimension not  exceeding
$n-c\delta^n$, where $c=c(n) >0$ is a  constant depending  only
on $n$. In particular a porous set has Lebesgue measure zero.

\vs 0,2cm\n Theorem 1.2 extends the same result established in \cite{[S]}
in the quasilinear degenerate and linear cases $p\geqslant 2$.
The proof is based on the exact growth of the solution of
the problem $(VF)$ near the free boundary which is given by the
next theorem. 
\vs 0.2cm 
\begin{theorem}\label{t1.3} Assume that $1<p<2$ and that $f$ and $g$ satisfy (1.1)-(1.4),
and let $u$ be a solution of the problem $(VP)$.  Then there exists two positive constants 
$c_0=c_0(n, p, \lambda_0)$ and $C_0=C_0(n, p, \lambda_0, \Lambda_0, \| g\|_{\infty,\Omega_T})$
such that for every compact set $K\subset\Omega_T$, $(x_0,t_0)\in (\partial\{u>0\})\cap K$,
the following estimates hold  
\begin{equation}\label{1.5}
c_0 r^q\leqslant \sup_{B_r(x_0)} u(.,t_0)\leqslant C_0 r^q,
\end{equation}
where $\displaystyle{q={p\over{p-1}}}$ is the conjugate of $p$. 
\end{theorem}

\vs 0,2cm\n Since the proof of Theorem 1.2 relies on the one of Theorem 1.3,
it will be enough to prove the latter one. On the other hand we observe that the left hand side 
inequality in (1.5) was established in \cite{[S]} Lemma 2.1 for any $p>1$, 
while the right hand side inequality in (1.5) was established only for 
$p\geqslant 2$. In the next section, we shall
establish the second inequality for a class of functions in the singular case 
i.e. for $1<p<2$. Then the right hand side inequality in (1.8) will follow 
exactly as in \cite{[S]} and we refer the reader to that reference for the details.
Hence the proof of Theorem 1.2 will follow.

\vs 0,3cm\n For similar results in the quasilinear elliptic case, we refer to
$\cite{[KKPS]}$, $\cite{[CL1]}$, and $\cite{[CL2]}$, respectively for the $p$-obstacle problem, 
the $A$-obstacle problem, and the $p(x)$-obstacle problem. For the obstacle problem for a
class of heterogeneous quasilinear elliptic operators with variable growth, we refer to \cite{[CLRT]}.

\section{A  class of functions on the unit cylinder}\label{2}

\vs 0.3cm \n In this section, we assume that $1<p<2$ and consider the family 
$\mathcal{F}=\mathcal{F}(p,n,M,\Lambda_0)$ of functions $u$ defined
on the unit cylinder $Q_1=B_1\times(-1,1)$ by
$u\in\mathcal{F}$ if it satisfies 
\begin{align}\label{e2.1-2.4}
   & u\in W^{1,p}(Q_1), \qquad  \|u_t-\Delta_p u \|_{L^\infty (Q_1)} \leqslant \Lambda_0\quad \hbox{ in } Q_1\\
   & 0\leqslant u \leqslant M  \quad \hbox{ in } Q_1 \\
   & u(0,0)=0\\
   &u_t\geqslant 0 \quad \hbox{ in } Q_1.
\end{align}

\vs 0,3cm\n The following theorem  gives the growth
of the  elements of the family $ \mathcal{F}$ in the singular case. 
This completes a result proved in \cite{[S]} for the degenerate case $p\geqslant 2$.

\begin{theorem}\label{t1.1} There exists a positive  constant
$C=C(p,n,M,\Lambda_0)$ such that for every $u\in \mathcal{F}$, we have
$$ u(x,t)\leqslant Cd(x,t) \qquad \forall (x,t)\in Q_{1/2}$$
where $d(x,t)=\sup\{r~/~Q_r(x,t)\subset \{u>0\}~\}$ for $(x,t)\in \{u>0\}$, and 
$d(x,t)=0$ otherwise, and where $Q_r(x,t)=B_r(y)\times(s-r^q,s+r^q)$.
\end{theorem}

\vs 0,2cm \n In order to prove Theorem 2.1, we need to introduce some
notations inspired from \cite{[S]}. For a nonnegative  bounded function  $u$, we define the
 quantities

$$Q_r^-=B_r\times(-r^q,0),\quad S(r,u) = \sup_{(x,t)\in Q_r^-} u(x,t).$$

\n We also define for $u\in\mathcal{F}$ the set
$$ \mathbb{M}(u) = \{ j\in \mathbb{N}\cup\{0\}/\quad A S( 2^{-j-1},u)
\geqslant S( 2^{-j},u) \}$$

\n where $A= 2^q\displaystyle{  \max\Big(1,{1\over {C_0}}\Big) } $  and
$C_0$ is the constant in (1.8). 

\n As in \cite{[S]}, we first show a weaker version of the inequality.

\vs 0.2cm \n
\begin{lemma}\label{l2.1} There exists  a constant $C_1=C_1(p,n,M,\Lambda_0)$ such
that
$$ S( 2^{-j-1},u) \leqslant C_1 2^{-qj} \qquad
 \forall u\in \mathcal{F} , \quad \forall j\in \mathbb{M}(u).$$
\end{lemma}

\vs 0.2cm\n \emph{Proof.} We argue by contradiction and assume that
\begin{equation}\label{e2.5}
\forall k\in  \mathbb{N},\quad \exists u_k \in \mathcal{F},\quad
  \exists j_k\in \mathbb{M}(u_k) \quad \hbox{ such that }\quad
  S( 2^{-j_k-1},u_k) \geqslant k 2^{-qj_k}.
\end{equation}
\n Let $\alpha_k=2^{-pj_k}(S(2^{-j_k-1},u_k))^{2-p}$,
and consider  $v_k(x,t)={{u_k(2^{-j_k}x,\alpha_k t)}\over{S(2^{-j_k-1},u_k)}}$ defined  in $Q_1$.  

\n First we observe that since $u(0,0)=0$ and $u$ is continuous, we have 
$\alpha_k\rightarrow 0$ as $k\rightarrow\infty$.

\vs 0.2cm\n Moreover, we have 
\begin{eqnarray}\label{e2.6-2.7}
&& \nabla v_k(x,t)={{2^{-j_k}}\over{S( 2^{-j_k-1},u_k)}}\nabla u_k(2^{-j_k}x,\alpha_k t)\nonumber\\
&& v_{kt}(x,t)={{\alpha_k}\over{S( 2^{-j_k-1},u_k)}}u_{kt}(2^{-j_k}x,\alpha_k t)
=\Big({{2^{-qj_k}}\over{S(2^{-j_k-1}},u_k)}\Big)^{p-1}u_{kt}(2^{-j_k}x,\alpha_k t)\\
&& \Delta_p v_k(x,t) =div\big(|\nabla v_k|^{p-2}\nabla v_k\big)\nonumber\\
&&\quad=\Big({{2^{-j_k}}\over{S(2^{-j_k-1}},u_k)}\Big)^{p-1}div\big(|\nabla u_k(2^{-j_k}x,\alpha_k t)|^{p-2}\nabla u_k(2^{-j_k}x,\alpha_k t)\big)\nonumber\\
&&\quad=2^{-j_k}\Big({{2^{-j_k}}\over{S(2^{-j_k-1}},u_k)}\Big)^{p-1}\Delta_p u_k(2^{-j_k}x,\alpha_k t)\nonumber\\
&&\quad=\Big({{2^{-qj_k}}\over{S(2^{-j_k-1}},u_k)}\Big)^{p-1}\Delta_p u_k(2^{-j_k}x,\alpha_k t).
\end{eqnarray}

\n We deduce from (2.6)-(2.7) that
\begin{equation}\label{e2.8}
v_{kt}-\Delta_p v_k(x,t) = \Big({{2^{-qj_k}}\over{S(2^{-j_k-1}},u_k)}\Big)^{p-1}(u_{kt}-\Delta_p u_k)(2^{-j_k}x,\alpha_k t).
\end{equation}

\n Combining (1.1), (2.1)-(2.5) and (2.8), we obtain

\begin{eqnarray}\label{e2.9-2.13}
&&  \|v_{kt}-\Delta_p v_k\|_{\infty} \leqslant {\Lambda_0\over{k^{p-1}}}\quad \hbox{ in } Q_1 \\
&& 0\leqslant v_k\leqslant{{S( 2^{-j_k},u_k)}\over{S( 2^{-j_k-1},u_k)}}\leqslant A\quad \hbox{ in } Q_1^-,\\
&&\qquad  v_{kt}\geqslant 0\quad \hbox{ in } Q_1^-,\\
&& \qquad \displaystyle {\sup_{Q_{1/2}^-} } v_k=1\\
&& \qquad v_k(0,t)=0\quad\forall t\in(-1,0).
\end{eqnarray}

\n Taking into account (2.9)-(2.10), we deduce (see \cite{[Ch]}) that there exists two
positive constants $\beta=\beta(n, p, M, A)$ and $C=C(n, p, M, A)$
such that $v_k\in C^{0,\beta}(\overline{Q}_{3/4}^-)\cap C_x^{1,\beta}(\overline{Q}_{3/4}^-)$ and
$$|v_k|_{\beta,\overline{Q^-_{3/4}}},~ |\nabla v_k|_{\beta,\overline{Q^-_{3/4}}} \leqslant C,\quad \forall k$$

\n It follows then from Ascoli-Arzella's theorem that there exists
a subsequence, still denoted by $v_k$ and a function $v\in
C^{0,\beta}(\overline{Q^-_{3/4}})\cap C_x^{1,\beta}(\overline{Q^-_{3/4}})$ such that $ v_k \longrightarrow
v$  and $ \nabla v_k \longrightarrow \nabla v$ uniformly in $\overline{Q^-_{3/4}}$. Moreover, using
(2.9)-(2.13), we see that $v$ satisfies

$$\left\{
    \begin{array}{ll}
 & v_t-\Delta_p v =0 \quad \hbox{ in } Q_{3/4}^-,\qquad  v, v_t\geqslant
0\quad \hbox{ in } Q_{3/4}^-,\\
&\\
&\displaystyle{\sup_{x\in Q_{1/2}^-} v(x,t)}=1,\qquad v(0,t)=0\quad\forall t\in(-3/4,0).
    \end{array}
  \right.
$$

\vs 0,3cm\n We discuss two cases:

\vs 0,3cm\n \emph{\underline{Case 1}:} $\forall (x,t)\in Q_{3/4}^-$  $v(x,t)=0$

\vs 0,3cm\n In particular we have $v\equiv 0$ in $Q_{1/2}^-$ which contradicts the fact that
$\displaystyle{\sup_{x\in Q_{1/2}^-} v(x)=1}$.

\vs 0,3cm\n \emph{\underline{Case 2}:} $\exists (x_0,t_0)\in Q_{3/4}^-$ such that $v(x_0,t_0)>0$

\vs 0,3cm\n Since $v(.,t_0)$ is not identically zero and $v(0,t_0/2)=0$, we get from the strong
maximum principle (see \cite{[N]}) that $v(x,t_0/2)=0$ for all $x\in B_{3/4}$. By the monotonicity
of $v$ with respect to $t$ and the fact that $v$ is nonnegative, we have necessarily $v(x,t)=0$ for all 
$(x,t)\in B_{3/4}\times(-3/4,t_0/2)$, which is in contradiction
with the fact that $v(x_0,t_0)>0$.

\qed

\vs 0,5cm \n \emph{Proof \, of \, Theorem \,2.1.} Using Lemma 2.1, the proof follows exactly as
the one of Theorem 2.2 in \cite{[S]}. 

\qed

\end{document}